\documentclass[oneside, a4paper,10pt]{amsart}
\usepackage{array,latexsym}
\usepackage[pdftex]{graphicx}
\usepackage{amsfonts}
\usepackage{amsmath}
\usepackage{amsthm}
\usepackage{fullpage}
\usepackage{enumerate}
\newtheorem{theorem}{Theorem}
\newtheorem{lemma}[theorem]{Lemma}

\newtheorem{corollary}[theorem]{Corollary}

\theoremstyle{definition}
\newtheorem{definition}{Definition}
\theoremstyle{remark}
\newtheorem{remark}{Remark}
\newtheorem{example}{Example}

\newcommand{\qh}[1]{[#1]_{q}}
\newcommand{\qxh}[1]{[#1]^{x}_{q}}

\newcommand{\pqxh}[1]{~_{p}[#1]^{x}_{q}}

\newcommand{\INV}[1]{\mathrm{inv}(#1)}
\newcommand{\NON}[1]{\mathrm{non}(#1)}
\newcommand{\CYC}[1]{\mathrm{cyc}(#1)}
\newcommand{\ASC}[1]{\mathrm{asc}(#1)}
\newcommand{\DES}[1]{\mathrm{des}(#1)}

\newcommand{\MAP}{\Theta}
\newcommand{\BIJ}{\Phi}
\newcommand{\PAR}{\Psi}

\newcommand{\ONE}{\Omega}
\newcommand{\DDD}{\Delta}

\newcommand{\opa}{r}
\newcommand{\opb}{\rho}
\newcommand{\opc}{\rho'}
\newcommand{\opd}{t}
\newcommand{\ope}{\tau}
\newcommand{\opf}{\tau'}

\begin{document}
\title{Cycles and patterns in permutations}
\author[R. Parviainen]{Robert Parviainen}
\address{ARC Centre of Excellence for Mathematics 
   	and Statistics of Complex Systems\\
  	139 Barry Street, The University of Melbourne, Victoria, 3010}
\email{robertp@ms.unimelb.edu.au}
\date{\today}
\maketitle
\begin{abstract}
We study joint distributions of cycles and patterns in permutations written in standard cycle form. Both non-restricted and generalised patterns, of length 2 and 3, are explored. Many extensions of classical theory are achieved; bivariate generating functions for inversions, ascents, descents, 123s, valleys, 1'-2-1s; closed forms for avoidance of peaks, 2-3-1s, 1-2-3s, 2'-1-2s and 1'-2-1s; bijective proofs of Wilf-equivalences.

We also derive some new results about standard pattern occurrence,  such as continued fractions for the generating functions for occurrences of valleys and for occurrences of the pattern 123.  

The methods are simple and combinatorial in nature: direct enumerative analysis and bijections to lattice paths.
\end{abstract}


\section{Introduction}
Not much is known about the interaction of cycle statistics and pattern occurrences in permutations. Some results are given by Edelman, \cite{Edel1987}, who studied cycles and inversions. Unlike Edelman, who considered inversions and cycles in permutations written in index form, we consider patterns in permutations in cycle form. This shift in viewpoint allows for a multitude of new results, which are often easily derived by simple methods.

Our results include
\begin{itemize}
	\item Continued fractions for generating functions for
	\begin{itemize}
		\item cycles and inversions
		\item cycles and ascents and descents
		\item cycles and double ascents
		\item cycles and valleys
	\end{itemize}
	\item Closed form expressions for
	\begin{itemize}
		\item avoidance of  peaks
		\item avoidance of 2--3--1
		\item avoidance of 1--2--3
		\item avoidance and 1 and 2 occurrences of 1-2-1'
		\item avoidance and 1 and 2 occurrences of 2-1-2'
	\end{itemize}
\end{itemize}

Many of the results are refinements of classical results: When considering occurrences of cycles and a certain pattern in cyclic form, the sum over the number of cycles is equivalent to the number of occurrences of the pattern in index form.

The rest of this section is devoted to general definitions and background. Sections 2 and 3 considers patterns of length two: inversions, and ascents and descents. Patterns of length three is the topic of sections 4, 5 and 6.

\subsection{Patterns in permutations}
Let $\mathcal{S}_{n}$ denote the set of permutations on $[n]=\{1,2,\dots, n\}$, and let $\mathcal S=\bigcup_{n}\mathcal{S}_{n}$. We write a permutation $\pi\in \mathcal{S}_{n}$ in one line index form as $\pi=\pi(1)\pi(2)\cdots\pi(n)$. A (generalised) \emph{pattern} is a permutation $\sigma\in\mathcal{S}_{k}$ and a set of restrictions. We first define occurrences of patterns without restrictions.

\begin{definition}
If $\sigma\in\mathcal{S}_{k}$ and $\pi\in \mathcal{S}_{n}$ we say that $\sigma$ occurs in $\pi$ if there exist $i_{1}<\dots < i_{k}$ such  that $\sigma=R(\pi(i_{1})\cdots\pi({i_{k}}))$, where $R$ is the reduction operator that maps the smallest element of the subword to 1, the second smallest to 2, and so on.	
\end{definition}

For example, an occurrence of the pattern 3--2--1 in $\pi\in \mathcal{S}_{n}$ means that there exist $1\leq i<j<k\leq n$ such that $\pi({i})>\pi({j})>\pi({k})$.

\emph{Generalised} or \emph {restricted} patterns were introduced by Babson and Steingr{\'\i}msson, \cite{BS2000}. The restriction is that two specified adjacent elements in the pattern \emph{must be adjacent} in the permutation as well. The position of the restriction in the pattern is indicated by \emph{an absence}  of a dash (--).  An occurrence of the pattern 3--21 in $\pi\in \mathcal{S}_{n}$ therefore means that there exist $1\leq i<j<n$ such that $\pi({i})>\pi({j})>\pi({j+1})$.

If a (generalised) pattern has no dashes the pattern is said to be a \emph{consecutive} pattern, as all elements must be consecutive in the permutation.

Introduced by Kitaev, see \cite{Kita2007}, \emph{partially ordered patterns} are a further generalisation. In a partially ordered pattern the letters form a partially ordered set. An occurrence of a partially ordered pattern is a linear extension of the corresponding partially ordered set in the order indicated by the pattern. Some simple partially ordered patterns are considered in sections \ref{section:212'}, \ref{section:121'}, \ref{section:2-1-2'} and \ref{section:1-2-1'}.

\subsection{Cycles}

\subsubsection{Standard cycle form}
The \emph{standard cycle form} of a permutation $\pi \in \mathcal{S}_{n}$ is the permutation written in cycle form, with cycles starting with the smallest element, and cycles ordered in decreasing order with respect to their minimal elements.  We write $(c_{1}\cdots c_{k})$ for a cycle with $k$ elements, and $\CYC{\pi}$ for the number of cycles in $\pi$. 
\begin{example}
We have $\pi=47613852=(275368)(14)$, and $\CYC{\pi}=2$.
\end{example}

The following standard bijections, denoted $\PAR$ and $\ONE$, between permutations in standard cycle form, and permutations and cycles and permutations, respectively, will occasionally be useful.
\begin{definition}
If the $\pi$ is given in standard cycle form as
\[
\pi=(c_{1}^{1} c_{2}^{1}\cdots c_{i_{1}}^{1})(c_{1}^{2} c_{2}^{2}\cdots c_{i_{2}}^{2})\cdots (c_{1}^{k} c_{2}^{k}\cdots c_{i_{k}}^{k}),
\]
define 
\[\PAR(\pi)=c_{1}^{1} c_{2}^{1}\cdots c_{i_{1}}^{1}c_{1}^{2} c_{2}^{2}\cdots c_{i_{2}}^{2}\cdots c_{1}^{k} c_{2}^{k}\cdots c_{i_{k}}^{k}\]
(i.e. the permutation achieved by removing the parenthesis).

If $c=(c_{1}c_{2}\cdots c_{k})$ is a cycle, define $\ONE(c)=c_{2}\cdots c_{k}$.
\end{definition}

\subsubsection{Cyclic occurrence of patterns}
\begin{definition}
Let $\pi$ be a permutation of $[n]$, with standard cycle form
\[
\pi=(c_{1}^{1} c_{2}^{1}\cdots c_{i_{1}}^{1})(c_{1}^{2} c_{2}^{2}\cdots c_{i_{2}}^{2})\cdots (c_{1}^{k} c_{2}^{k}\cdots c_{i_{k}}^{k}),
\]
and let $\sigma$ be a generalised pattern. It occurs \emph{cyclically} in $\pi$ if it occurs in the permutation
\[
\PAR(\pi)=c_{1}^{1} c_{2}^{1}\cdots c_{i_{1}}^{1}c_{1}^{2} c_{2}^{2}\cdots c_{i_{2}}^{2}\cdots c_{1}^{k} c_{2}^{k}\cdots c_{i_{k}}^{k},
\]
with the further restriction that $c_{i_{j}}^{j}$ and $c_{1}^{j+1}$ are \emph{not} adjacent. For example, 21 does not occur in (2)(1), but 2--1 does.
\end{definition}

\subsection{Motzkin paths}
\begin{definition}
  A \emph{Motzkin path} of length $n$ is a sequence of  vertices $p = (v_0,v_1,\dots,v_n)$, with $v_i \in \mathbb{N}^2$  (where $\mathbb{N} = \{0,1,\dots\}$), with steps $v_{i+1}-v_i \in \{ (1,1), (1,-1), (1,0)\}$ and $v_0 = (0,0)$ and $v_n=(n,0)$. 
  
   A  \emph{bicoloured Motzkin path} is a Motzkin path in which each east, $(1,0)$, 
  step is labelled by one of two colours. 
\end{definition}

\emph{Unless otherwise stated, all Motzkin paths studied here will be bi-coloured}, and we denote the set of paths of length $n$ by $\mathcal{M}_{n}$. 

\begin{remark}
Although the use of bi-coloured paths is not necessary, it simplifies some discussions. 
\end{remark}

Let $\mathrm{N}$ ($\mathrm{S}$) denote a north, $(1,1)$, step (resp., south, $(1,-1)$, step), and $\mathrm{E}$ and $\mathrm{F}$ the two differently coloured  east steps. Further, let $N_{h}, S_{h}, E_{h}, F_{h}$ denote the weights of  $\mathrm{N}$, $\mathrm{S}$, $\mathrm{E}$, $\mathrm{F}$ steps, respectively, that start at height $h$. The weight of a Motzkin path is the product of the steps' weights. 

Very useful is Flajolet's \cite[Theorem 1]{Flaj1980} continued fraction representation for Motzkin path generating functions.

\begin{theorem}
Let $w(p)$ denote the weight of a Motzkin path $p$. Then
\begin{multline*}
\sum_{n}\sum_{p\in\mathcal{M}_{n}}w(p)z^{n}\\
=\cfrac{1}{1-(E_{0}+F_{0})z-
\cfrac{N_{0}S_{1}z^{2}}{1-(E_{1}+F_{1})z-
\cfrac{N_{1}S_{2}z^{2}}{1-(E_{2}+F_{2})z-
\cfrac{N_{2}S_{3}z^{2}}{1-(E_{3}+F_{3})z-\cdots}}}}
\end{multline*}
\end{theorem}

\subsection{Further notation}
If $f_{p}(k,m,n)$ denotes the number of permutations of length $n$ with $m$ cycles and $k$ occurrences of a pattern $p$, 
the \emph{(ordinary) generating function} for cyclic occurrences of the pattern $p$ is 
\[F_{p}(q,x,z)=\sum_{m,k,n}f_{p}(k,m,n)q^{k}x^{m}z^{n}.\]

As in the example above, we will use $z$ to mark permutation length, $x$ to mark number of cycles, and $q$ to mark occurrences of some pattern, and $f_{p}(k,m,n)$ and $F_{p}(q,x,z)$ for the counts and generating functions, respectively. 

For the case of (classic) non-cyclic occurrences, let $g_{p}(k,n)$ denote the number of permutations of length $n$ with $k$ occurrences of the pattern $p$, and let $G(q,z)=\sum_{k,n}g_{p}(q,z)$  denote the generating function.

Finally, we define the $q$-number $\qh{h}$ and some generalisations. Let
\begin{align*}
\qh{h}&=1+q+\cdots + q^{h-2} + q^{h-1},\\
\qxh{h}&=1+q+\cdots + q^{h-2}+ {x} q^{h-1},\mbox{ and}\\
\pqxh{h}&=p^{h-1}+p^{h-2}q+\cdots + p q^{h-2}+ {x} q^{h-1}.
\end{align*}


\section{Inversions}
An inversion in a permutation is a pair  $i<j$ such that $\pi_{i}>\pi_{j}$. In terms of patterns, it is an occurrence of the pattern 2--1. Let $\INV{\pi}$ denote the number of cyclic inversions in a permutation $\pi$, and $\NON{\pi}=\binom{|\pi|}{2}-\INV{\pi}$ the number of non-inversions.

It is well known, \cite{Muir1898}, that the generating function $G_{\mbox{\scriptsize 2--1}}(q,z)=\sum_{\pi\in\mathcal S}q^{\INV{\pi}} z^{|\pi |}$ for the number of permutations with $k$ inversions is 
\[G_{\mbox{\scriptsize 2--1}}(q,z)=\sum_{n=1}^{\infty}z^{n }\prod_{i=1}^{n}\qh{i}.\]

A straightforward extension gives the 4-variable generating function
\[F_{\mbox{\scriptsize 1--2, 2--1}}(p,q,x,z)=\sum_{\pi\in\mathcal{S}}p^{\NON{\pi}}q^{\INV{\pi}}x^{\CYC{\pi}}z^{|\pi |}.\]
\begin{theorem}
\begin{align*}
F_{\mbox{\scriptsize 1--2, 2--1}}(p,q,x,z)=&\sum_{n=1}^{\infty}z^{n }\prod_{i=1}^{n}\pqxh{i}\\
		=&
		\cfrac{\pqxh{1}z}{1-
		\cfrac{\pqxh{2}z}{1+\pqxh{2}z-
		\cfrac{\pqxh{3}z}{1+\pqxh{3}z-
		\cfrac{\pqxh{4}z}{1+\pqxh{4}z+\cdots}}}}.
\end{align*}
\end{theorem}

\begin{proof}
Induction on permutation length. Let $F_{n}(p,q,x)=\sum_{\pi\in\mathcal{S}_{n}}p^{\NON{\pi}}q^{\INV{\pi}}x^{\CYC{\pi}}$. It is trivial to check that $F_{1}(p,q,x)=1$.

Assume that $F_{n}(p,q,z)=\prod_{i=1}^{n}\pqxh{i}$, and consider $\pi\in \mathcal{S}_{n}$. Write the permutation in cycle form:
\[C(\pi)=(c_{1}^{1} c_{2}^{1}\cdots c_{i_{1}}^{1})(c_{1}^{2} c_{2}^{2}\cdots c_{i_{2}}^{2})\cdots (c_{1}^{k} c_{2}^{k}\cdots c_{i_{k}}^{k}).\]
Consider the insertion of $n+1$. We can either add the cycle $(n+1)$, which adds $n$ inversions and one cycle, or insert $n+1$ into one of the cycles. There are $n$ valid choices; after each element in each cycle (we cannot put $n+1$ first in a cycle as all cycles starts with their least element). Inserting after $c_{1}^{1}$ adds $n-1$ inversions and one non-inversion, after $c_{2}^{1}$ adds $n-2$ inversions and 2 non-inversions, and so on. 

Summing up, $F_{n+1}(p,q,x)=(p^{n}+p^{n-1}q+\cdots+pq^{n-1}+xq^{n})F_{n}(p,q,x)=\prod_{i=1}^{n+1}\pqxh{i}$.

The second equality is an immediate consequence of a result of Euler \cite{Eule1775}.
\end{proof}


\section{Ascents and descents, or 12 and 21}
An ascent in a permutation is an index $i$ such that $\pi_{i}<\pi_{i+1}$. In terms of patterns, it is an occurrence of the pattern 12. A descent is an index $i$ such that $\pi_{i}>\pi_{i+1}$ --- it is an occurrence of the pattern 21. Let $\ASC{\pi}$ and $\DES{\pi}$ denote the number of ascents and descents, respectively, in a permutation $\pi$.

The next result is a refinement of the standard continued fraction representation of the Eulerian numbers.
\begin{theorem}\label{theorem:asc+des}
\begin{multline*}
F_{\mbox{\scriptsize 21, 12}}(p,q,x,z)=
\sum_{\pi\in \mathcal S}p^{\DES{\pi}}q^{\ASC{\pi}}x^{\CYC{\pi}}z^{|\pi |}\\
=\cfrac{1}{1-x z-
\cfrac{q x z^{2}}{1-(q+p+x)z-
\cfrac{2 q (p+x) z^{2}}{1-(2q+2p+x)z-
\cfrac{3 q (2p+x) z^{2}}{1-(3q+3p+x)z-\cdots}}}}
\end{multline*}
\end{theorem}

\subsection{Proof of Theorem \ref{theorem:asc+des}}
We want to use the fact that the continued fraction is the generating function for weighted Motzkin paths with weights
\[N_{h}=(h+1)q, S_{h}=(h-1)p +x, E_{h}=hq+x, F_{h}=h p.\]
To this end we will give a bijection between permutations and Motzkin paths.

\subsection{The arc diagram representation}
\begin{figure}[htbp]
\begin{center}
\includegraphics[scale=0.8]{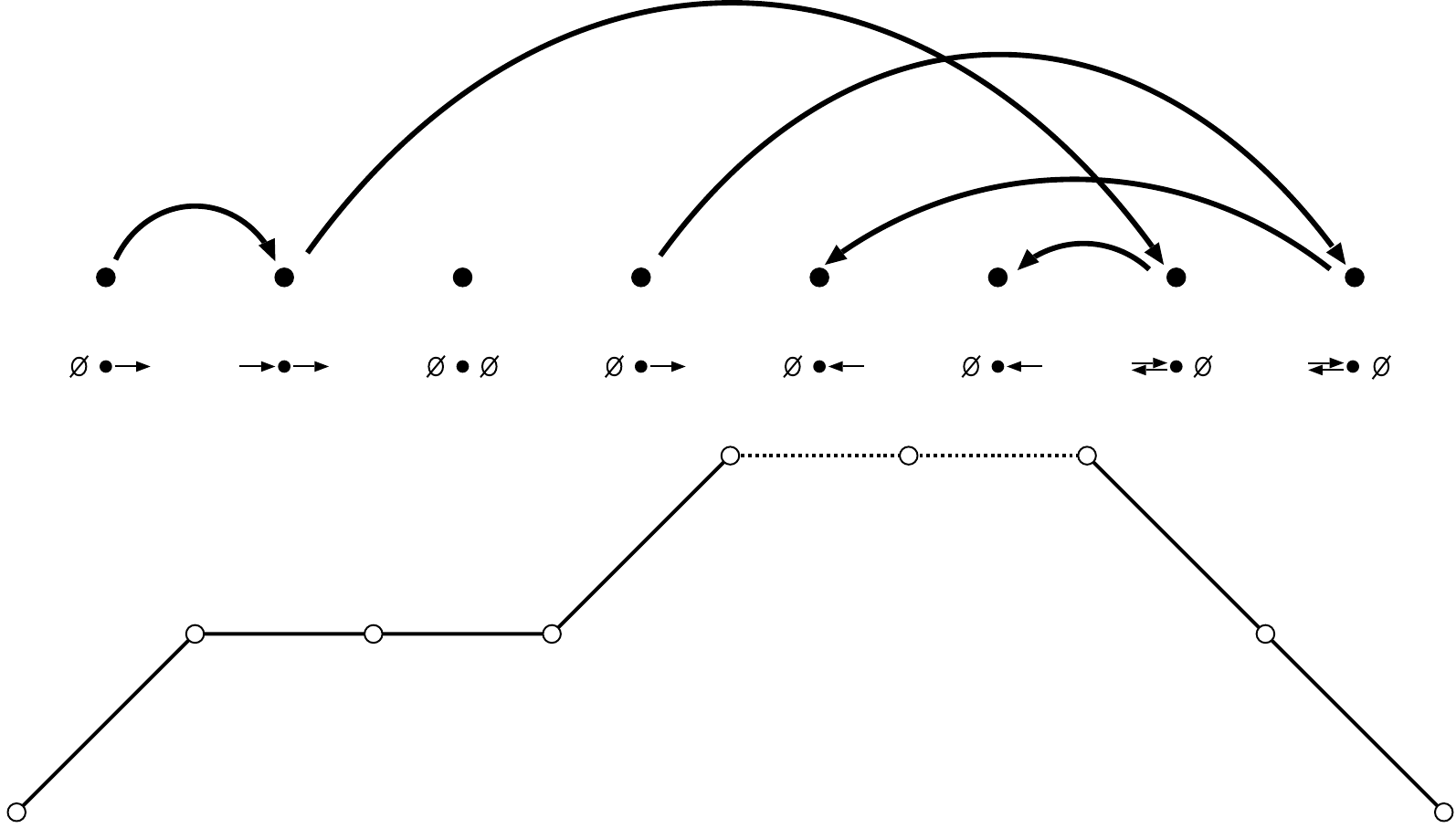}
\caption{The arc diagram for the permutation $\pi=(485)(3)(1276)$, the shape pairs, and the path $\MAP(\pi)$ ($\mathrm{E}$ steps are solid, $\mathrm{F}$ steps dotted).}
\label{fig:1}
\end{center}
\end{figure}

We use a graphical representation of permutations to aid in the description of the mapping. For permutation $\pi\in\mathcal{S}_n$ with standard cycle form 
\[
\pi=(c_{1}^{1} c_{2}^{1}\cdots c_{i_{1}}^{1})(c_{1}^{2} c_{2}^{2}\cdots c_{i_{2}}^{2})\cdots (c_{1}^{k} c_{2}^{k}\cdots c_{i_{k}}^{k}),
\]
put $n$ nodes on a line, representing the elements $1$ to $n$. For $s=1,\ldots, k$ and $t=1,\ldots, i_{s}-1$ draw an arc from node $c_{t}^{s}$ to node $c_{t+1}^{s}$.  See Figure \ref{fig:1}  for an example.

Associate each node with a left and a right shape. The left (right) shape is the set of connections to nodes on the left (right) side with the node. The possible shapes on both sides are $\{\emptyset, \rightarrow, \leftarrow, \rightleftharpoons\}$. See Figure \ref{fig:1} again. 
Give node $k$ weight $x$ if it is the rightmost in a cycle, weight $q$ if there is an arc leaving the node to the right, and weight $p$ otherwise. 

We now define a mapping from arc diagrams to paths, which will be the base for our bijection.

\begin{definition}
If $\pi\in\mathcal{S}_{n}$ have left shapes $\{l_{1},\ldots, l_{n}\}$ and right shapes $\{r_{1},\ldots, r_{n}\}$, let step $k$ in $\MAP(\pi)$ be $s_{k}$, where $s_{k}$ is given by the following table (where ``$-$'' denotes shape pairs that do not appear). Further, give step $k$ the same weight as node $k$.
\begin{center}
	\begin{tabular}{c||c|c|c|c}
		$l_{k}\backslash r_{k}$ & $\emptyset$ & $\rightarrow$ & $\leftarrow$ & $\rightleftharpoons$ \\\hline\hline 
		$\emptyset$ & $\mathrm{E}$ & $\mathrm{N}$ & $\mathrm{F}$ & $\mathrm{N}$ \\\hline 
		$\rightarrow$ & $\mathrm{S}$ & $\mathrm{E}$ & $-$ & $-$ \\\hline 
		$\leftarrow$ & $-$ & $-$ & $\mathrm{F}$ & $-$ \\\hline 
		$\rightleftharpoons$ & $\mathrm{S}$ & $-$ & $-$ & $-$
	\end{tabular}
\end{center}
\end{definition}

\begin{lemma}
The mapping  $\MAP$ is a surjection from the set of permutations to the set of Motzkin paths with no $\mathrm{F}$ steps at level 0.
\end{lemma}
\begin{proof}
To show that the image is a Motzkin path, the conditions for a Motzkin path must be verified. Namely, that $N_{k}\geq S_{k}$, $k<n$, and $N_{n}=S_{n}$, where $N_{m}$ and $S_{m}$ are the number of $\mathrm{N}$ and $\mathrm{S}$ steps, respectively, up to and including step $m$.

As the shape pairs $(\leftarrow, \leftarrow)$,  $(\rightarrow, \rightarrow)$, $(\emptyset, \emptyset)$ and $(\emptyset, \leftarrow)$ map to $\mathrm{E}$ and $\mathrm{F}$ steps, it may be assumed that these shapes do not occur.

Now, in a valid arc diagram, the number of $(\emptyset, \rightarrow)$ and $(\emptyset, \rightleftharpoons)$ shape pairs up to and including node $k$ must be greater than or equal to the number of $(\rightarrow, \emptyset)$ and $(\rightleftharpoons, \emptyset)$ shape pairs, because the incoming arcs in the latter two shape pairs must start somewhere. Further, these counts must agree for $k=n$. This is exactly what is needed. 

To show that $\MAP$ is a surjection, consider any Motzkin path $r$ with no $\mathrm{F}$ steps at level 0. We will build an arc diagram $a$ that maps to $r$. 

For each $\mathrm{F}$ step in $p$ we can associate a unique pair of $\mathrm{N}$ and $\mathrm{S}$ steps. (The rightmost (leftmost) $\mathrm{N}$ ($\mathrm{S}$) step to the left (right) of the $\mathrm{F}$ step, ending (starting) at the $\mathrm{F}$ steps level.) Let $n, f, s$ denote the positions of the $N, F, S$ steps, respectively. In the arc diagram $a$, draw an arc from node $n$ to node $s$ and one from node $s$ to node $f$. 

For the remaining $\mathrm{N}$ and $\mathrm{S}$ steps, fix one pairing of these, and draw arcs from the nodes corresponding to the $\mathrm{N}$ steps, to the associated nodes corresponding to the $\mathrm{S}$ steps.
 
Clearly, $a$ represents a permutation, and $\MAP(a)=r$ as desired.
\end{proof}

The next step is to show that $\MAP$ defines a bijection $\BIJ$ between the set of equivalence classes of permutations and the set of weighted Motzkin paths, where two permutations are equivalent if they map to the same \emph{unweighted} Motzkin path.
\begin{definition}
	 For an equivalence class $\mathbf{E}_{r}=\{\pi\in \mathcal{S}| \MAP(\pi)=r\}$ of permutations let $\BIJ(\mathbf{E}_{r})=r$, and let the weight of step $k$ be the sum of weights of node $k$ over permutations in $\mathbf{E_{r}}$.
\end{definition}
\begin{theorem}
The mapping $\BIJ$ is a bijection from the set of equivalence classes of permutations (with the above definition of equivalent) to the set of weighted Motzkin paths, with weights
\begin{equation*}\label{eq:weights 2} 
N_{h}=(h+1)q, S_{h}=(h-1) p +x, E_{h}=hq+x, F_{h}=h p.
\end{equation*}
such that the sum of weights of permutations in an equivalence class $\mathbf{E}$ is the weight of $\BIJ(\mathbf{E})$.
\end{theorem}
\begin{proof}
	That $\BIJ$ is a bijection follows at once from the fact that $\MAP$ is into the set of Motzkin paths, and $\BIJ$ is defined from the set of equivalence classes that maps to the same Motzkin paths.

	Fix a Motzkin path $\gamma$ and let $\mathbf{E}$ be an equivalence class of permutations that maps to it. Each permutation in $\mathbf{E}$ has a common \emph{skeleton}, the sequence of node types: if step $k$ in $\gamma$ is $N$, node $k$ in the arc diagram of any $\pi\in\mathbf{E}$ is of type $(\emptyset, \rightarrow)$ or $(\emptyset, \rightleftharpoons)$, etc. 
	
	We must show that with weights given by \eqref{eq:weights 2}, the weight of the path is given by the sum of weights of the permutations in $E$. We will do this by a weighted count of the number of ways the skeleton can be completed by filling in the arcs.	
	
	Suppose that the path has no horizontal ($\mathrm{E}$ or $\mathrm{F}$) steps. 
	
	First we connect arcs going from right to left. Consider the $S$ steps in order, and let the heights of the steps be $h_1, h_2, \ldots, h_k$. To the left of the first $S$ step there are $h_1$ $N$ steps. The node corresponding to the first step is either of type $(\rightarrow, \emptyset)$ or $(\rightleftharpoons, \emptyset)$. In the latter case there is $h_1-1$ possibilities for the endpoint of the arc originating at the node (the leftmost node cannot be of type $(\emptyset, \rightleftharpoons)$). The total weight of the possibilities is $(h_1-1)p+x$. Considering the next $S$ step, we find similarly that the total weight of the possibilities is $(h_2-1)p+x$ (remember that there already is one arc filled in).
	
	To fill in the arcs going from left to right, we go over the $N$ steps in the same way, but this time starting at the right end. Now it is found that the contribution of an $N$ step at height $h$ is $(h+1)q$.
	
We still have to consider $\mathrm{E}$ and $\mathrm{F}$ steps. Assume we have a fully connected arc diagram, and add an $\mathrm{E}$ step at height $h$. It may be seen that there are exactly $h$ arcs passing over the node from right to left. Thus, there are $h+1$ possible ways of connecting the node --- by adding it into one of the $h$ arcs (each give a node of type $(\rightarrow, \rightarrow)$ with weight $q$) or as an isolated node with weight $x$. So an $\mathrm{E}$ step at height $h$ contributes $hq+x$ towards the total. In the same way it can be shown that the contribution of $\mathrm{F}$ steps at height $h$ is $hp$.
\end{proof}


\section{Consecutive patterns}
In an occurrence of a consecutive pattern the letters have to be consecutive. Earlier we studied 12 and 21, here we look at consecutive patterns of length 3.
\subsection{123}
The series starts like this:
\begin{align*}
F_{123}(q,x,z)=&1 + x z + 
(x + x^2) z^2 + 
((1 + q) x + 3 x^2 + x^3) z^3\\ 
&+\left((3 + 2 q + q^2) x + (7 + 4 q) x^2 + 6 x^3 + x^4\right)z^4\\
&+ \big((9+11q+3q^2+q^3)x+(25+20q+5q^2)x^{2}\\
&\quad +(25+10q)x^{3}+10x^4+x^5\big)z^5+O(z^{6}).
\end{align*}
The approach used for ascents and descents again gives a lattice path interpretation, which then gives a continued fraction representation of the generating function.
\begin{theorem}\label{theorem:123}
The function $F_{123}(q,x,z)$ is the generating function for Motzkin paths with weights
\begin{align*}
E_{h}+F_{h}&=h (1 + q) + x\mbox{ and}\\
N_{h}S_{h+1}&=(h + 1)(h + x).
\end{align*}
\end{theorem}
\begin{proof}
We want to use bijection $\BIJ$ again. Hence, give node $k$ weight $x$ if it  is the rightmost in a cycle, weight $q$ if it is of the type $(\rightarrow, \rightarrow)$, and weight 1 otherwise. The proof now proceeds as the proof of Theorem \ref{theorem:asc+des}, and is omitted.
\end{proof}
\begin{corollary}
\begin{equation}
	F_{123}(q,x,t)=
	\cfrac{1}{1-xz
 	-\cfrac{xz^{2}}{1-(1+q+x)z
	-\cfrac{2(1+x)z^{2}}{1-(2+2q+x)z
 	-\cfrac{3(2+x)z^{2}}{1-(3+3q+x)z
 	\dotsb}}}}.
\end{equation}
\end{corollary}
As a further corollary we get the (ordinary) generating function for (non-cyclic) occurrences of 123 --- the exponential generating function were given in \cite{EN2003}.
\begin{corollary}
The ordinary generating function for occurrences of 123 is given by $F_{123}(q,1,z)$.
\end{corollary}
\begin{proof}
Follows by applying bijection $\PAR$. Cyclic occurrences of 123 in a permutation $\pi$ is in a one-to-one correspondence with with non-cyclic occurrences in  $\PAR(\pi)$.
\end{proof}

\subsection{132}
We will have little to say about 132, except giving the first few terms in the series,
\begin{align*}
F_{132}(q,x,z)=
&1 + x z + 
(x + x^2) z^2 + 
\left((1 + q)x + 3x^{2} + x^{3}\right)z^{3}\\
&+\left((2 + 4q)x + (7 + 4q)x^{2} + 6x^3 + x^4\right)z^{4}\\
&+\left((7+14q+3q^{2})x+(20+30q)x^{2}+25x^{3}+10qx^{3}
	+10x^{4}+x^{5}\right)z^{5}+O(z^{6}).
\end{align*}

One might ask if there is a simple continued fraction, similar to those seen before, that generates this series. We have not been able to find such a representation. (The same holds for all cases below where no such continued fraction representation is given.) By equating coefficients of $z$ one can solve for unknowns in the equation
\begin{equation*}
	F_{p}(q,x,t)=
	\cfrac{1}{1-a_0z
 	-\cfrac{b_0z^{2}}{1-a_1 z
	-\cfrac{b_1z^{2}}{1-a_2 z
 	-\cfrac{b_2z^{2}}{1-a_3 z
 	\dotsb}}}},
\end{equation*}
and hope for simple polynomial solutions from which a conjecture may be inferred. In the unsolved cases, it appears that the solutions are either rational functions or complicated polynomials including negative terms, and any simple pattern has eluded the author.

\subsection{213, 231, and 312}
The common series start 
\begin{align*}
F_{213}(q,x,z)=&F_{231}(q,x,z)=F_{312}(q,x,z)=
1 + x z + 
(x + x^2) z^2 + 
(2x + 3x^{2} + x^{3})z^{3}\\
&+\left((5 + q)x + 11x^2 + 6x^3 + x^4\right)z^{4}\\
&+\left((16+8q)x+(45+5q)x^{2}+35x^{3}+10x^{4}+x^{5}\right)z^{5}+O(z^{6}).
\end{align*}
Indeed, the patterns are equidistributed.
\begin{theorem}\label{theorem:213}
\[F_{213}(q,x,z)=F_{231}(q,x,z)=F_{312}(q,x,z).\]
\end{theorem}
\begin{proof}
Since the patterns are consecutive, it suffices to look at a single cycle. But cycles start with their lowest element, which can never be part of a 213, 231 or 312 pattern. The result therefore follows from the corresponding non-cyclic statement (which follows at once using the elementary permutation mappings reflection and complement), and an application of bijection $\ONE$.
\end{proof}
From the reasoning in the above proof it also follows that for these three patterns, non-cyclic occurrences are equivalent to cyclic occurrences in permutations with one cycle.
\begin{corollary}\label{corollary:213}
For $p\in\{213, 231, 312\} $,
\[f_{p}(k,1,n)=g_{p}(k,n-1).\]
\end{corollary}

\subsection{321} 
As for 132, we will not have much to say.
\begin{align*}
F_{321}(q,x,z)=
&1 + x z + 
(x + x^2) z^2 + 
(2x + 3x^{2} + x^{3})z^{3}\\
&+\left((5 + q)x + 11x^2 + 6x^3 + x^4\right)z^{4}\\
&+\left((17+6q+q^{2})x+(45+5q)x^{2}+35x^{3}+10x^{4}+x^{5}\right)z^{5}+O(z^{6}).
\end{align*}
Again, non-cyclic occurrences are equivalent to cyclic occurrences in permutations with one cycle.
\begin{corollary}
\[f_{321}(k,1,n)=g_{321}(k,n-1).\]
\end{corollary}

\subsection{Valleys, or 212', or ``213 or 312''}\label{section:212'}
Let $k=k_{1}+k_{2}$ and let $f_{213,312}(k,m,n)$ denote the number of permutations in $\mathcal S_{n}$ with $m$ cycles and  $k_{1}$ occurrences of 213 and $k_{2}$ occurrences of 312. 

\begin{remark}
This is the same as considering occurrences of the partially ordered pattern 212', \cite{Kita2007}. I.e.~ there exists $i<j<k$ such that $\pi(i)>\pi(j)<\pi(k)$, but there is no restriction on the relative magnitude of $\pi(i)$ versus $\pi(k)$.
\end{remark}

The series starts
\begin{align*}
F_{213,312}(q,x,z)=
&1 + x z + 
(x + x^2) z^2 + 
(2x + 3x^{2} + x^{3})z^{3}\\
&+\left((4 + 2q)x + 11x^2 + 6x^3 + x^4\right)z^{4}\\
&+\left((8+16q)x+(40+10q)x^{2}+35x^{3}+10x^{4}+x^{5}\right)z^{5}+O(z^{6}).
\end{align*}

Once again using the idea for ascents, a lattice path interpretation, and a continued fraction for the generating function can be found.
\begin{theorem}
The function $F_{213,312}(q,x,z)$ is the generating function for Motzkin paths with weights
\begin{align*}
E_{h}+F_{h}&=2h + x\mbox{ and}\\
N_{h}S_{h+1}&=(h + 1)(h q + x).
\end{align*}
\end{theorem}
\begin{proof}
The proof is similar to the proofs of Theorems \ref{theorem:asc+des} and \ref{theorem:123}, and omitted.
\end{proof}

\begin{corollary}
\begin{equation*}
	F_{213,312}(q,x,t)=
	\cfrac{1}{1-xz
 	-\cfrac{xz^{2}}{1-(2+x)
	-\cfrac{2(q+x)z^{2}}{1-(4+x)
 	-\cfrac{3(2q+x)z^{2}}{1-(6+x)
 	\dotsb}}}}.
\end{equation*}

\end{corollary}
We also get the (ordinary) generating function for (non-cyclic) occurrences of valleys. A recursion for the generating functions for $k$ occurrences were given by Rieper and Zekele, \cite{RZ2000}. Here we get a continued fraction for the full bivariate generating function. 
\begin{corollary}
The ordinary generating function for occurrences of valleys is given by 
\[G_{213, 312}(q,z)=\frac{F_{213, 312}(q,x,z)-1}{xz}\bigg| _{x=0}.\]
\end{corollary}
The above result will follow after the next lemma --- another simple consequence of the bijection $\ONE$ between $n+1$-cycles and $n$-permutations.
\begin{lemma}
The number of permutations of length $n$ with $k$ valleys equals the number of cycles of length $(n+1)$ with $k$ valleys. 
\end{lemma}

\subsection{Peaks, or 121', or ``132 or 231''}\label{section:121'}
As above we consider occurrences of a partially ordered pattern, 121'. This is equivalent to studying occurrences of both 132 and 231. The series starts
\begin{align*}
F_{132,231}(q,x,z)=
&1 + x z + 
(x + x^2) z^2 + 
\left((1 + q)x + 3x^{2} + x^{3}\right)z^{3}\\
&+\left((1 + 5q)x + (7 + 4q)x^{2} + 6x^3 + x^4\right)z^{4}\\
&+\left((x+18q+5q^{2})x+(15+35q)x^{2}+25x^{3}+10qx^{3}
	+10x^{4}+x^{5}\right)z^{5}+O(z^{6}).
\end{align*}
\begin{theorem}
Let $S(n,k)$ denote the Stirling numbers of the second kind, $S(n,k)=$ the number of partions of $[n]$ into $k$ non-empty parts. Then
\[f_{132,231}(0,k,n)=S(n, k).\]
\end{theorem}
\begin{proof}
Since cycles start with their lowest element, if a cycle avoids 132 and 231 it must be increasing. Therefore the number of permutations with $k$ cycles that avoids 132 and 231 equals the number of partions of $[n]$ into $k$ non-empty parts.
\end{proof}


\section{Patterns with one dash}
We refer to \cite{CM2002} for a summary of result conserning non-cyclic occurrences, and \cite{Parv2006b} for some results for cyclic occurrences. We just mention that in one case, that of 2--13, we have a continued fraction form for the generating function.
\begin{theorem}[\cite{Parv2006b}]
\begin{align*}
	F_{\mbox{\scriptsize 2--13}}(q,x,t)=
	\cfrac{1}{1-t (\qh{1}+\qxh{1})
   	-\cfrac{t^{2}\qh{1}\qxh{2}}{1-t (\qh{2}+\qxh{2})
	-\cfrac{t^{2}\qh{2}\qxh{3}}{1-t (\qh{3}+\qxh{3}) \dotsb}}}.
\end{align*}
\end{theorem}


\section{Unrestricted patterns}
For patterns having no restrictions, cyclic occurrence is a direct refinement of non-cyclic occurrence, in the following sense.
\begin{theorem}\label{theorem:nodash}
If a pattern $p$ is unrestriced, then 
\[\sum_{m=1}^{n}f_{p}(k,m,n)=g_{p}(k,n).\]
\end{theorem}
\begin{proof}
Follows by applying bijection $\PAR$; occurrences of unrestricted patterns in the standard cycle form of a permutation $\pi$ is in a one-to-one correspondance with occurrences in $\PAR(\pi)$.
\end{proof}

\subsection{Some bijections}
The aim of this section is to bijectively  demonstrate these claims.
\begin{enumerate}[\it i.]
	\item The distribution of cycles in 1--3--2-avoiding and 1--2--3-avoiding permutations are given by the Narayana numbers.
	\item The distribution of cycles in 2--1--3-avoiding and 3--1--2-avoiding permutations are given by the Catalan triangle. 
	\item The the distribution of cycles and occurrences of 3--1--2 and that of cycles and occurrences of 2--3--1 are equivalent. 
\end{enumerate}
Claim \emph{iii} is proved by a suitable involution. We also wish to prove \emph{i} and \emph{ii} bijectively. There are several bijections between 1--3--2-avoiding and 1--2--3-avoiding permutations in the literature \cite{EN2004, Krat2001, MDD2006, Reif2003, SS1985, West1995}, but none of them fully preserve the cycle structure (the bijection in \cite{EN2004} preserves the number of fixed points). We implicitly define bijections below, between  1--3--2-avoiding and 1--2--3-avoiding permutations, via Dyck paths, and between 2--1--3-avoiding and 3--1--2-avoiding permutations, again via Dyck paths. These have the nice property that the number of cycles maps to peaks and returns in the Dyck paths. They can be seen to be essentially variants of a two part bijection given by Krattenthaler, \cite{Krat2001}.

A \emph{Dyck path} is simply a Motzkin path with no $\mathrm{E}$ and $\mathrm{F}$ steps. A decent is any maximal substring of consecutive $S$ steps. A peak is an occurrence of an $N$ step immediately followed by a $S$ step. A return is  any visit to the $y$-axis, except the first. An excursion is a subwalk between to consecutive visits to the $y$-axis.

\subsubsection{Equivalence of avoidance of 1--2--3 and 1--3--2}
Three sets of operators is to defined. Each operator will construct a Dyck path or a permutation of size $n+1$ from one of size $n$.

\subsubsection*{Operator set 1 for Dyck paths}
If a path $p$ decomposes into excursions as $p=e_1e_2\cdots e_k$, the operators $\opa_i$, $i=1,\ldots, k$ are defined as $\opa_i p=e_1\cdots e_{i-1} N e_{i} \cdots e_k S$. That is, leave the first $i-1$ excursions unchanged, and raise the rest of the path. Also define an operator $\opa_0$ by $\opa_0 p=p NS$. 

\subsubsection*{Operator set for 1--2--3}
Assume a permutation $\pi$ of $[n]$ has cycle form $\pi=c_1c_2\cdots c_k\cdots c_m$, where $c_k$ is first cycle of size 2 or larger. The operators $\opb_i$, $i=1,\ldots, k$ are defined as $\opb_i \pi=c_1 \cdots c_{i-1} \hat c_i c_{i+1}\cdots c_m$, where $\hat c_i$ is the cycle $c_i$ with $n+1$ inserted immediately after the first (smallest) element in $c_i$. Also define $\opb_0$ by $\opb_0 \pi=(n+1)\pi$.

\begin{lemma}\label{lem:generate123}
	The set of permutations generated by $\opb_i$ are exactly those which avoid 1--2--3.
\end{lemma}
\begin{proof}\label{pf:generate123}
 None of the operators create an increasing subsequence of length 3. Conversely, if a permutation avoids 1--2--3, it can be built up by successive applications of $\opb_i$: Look at the element $n$ in a permutation $\pi$ of $[n]$. Since it is 1--2--3-avoiding, $n$ must either be a 1-cycle, or immediately after the smallest element in a cycle of length 2 or longer. Therefore the exists an $i$ such that $\pi=\opb_i \tilde \pi$ for some $\tilde \pi\in\mathcal{S}_{n-1}$.
\end{proof}

\subsubsection*{Operator set for 1--3--2}
Assume a permutation $\pi$ of $[n]$ has cycle form $\pi=c_1\cdots c_m$. Let $i_1,\ldots, i_k$ be the indices of the cycles whose smallest element is greater than all elements in all cycles to the right. The operators $\opc_i$, $i=1,\ldots, k$, inserts $n+1$ last in cycle $i_{k+1-i}$. Also, let $\opc_0$ be defined as $\opc_0 \pi=(n+1)\pi$. 

\begin{lemma}
	The set of permutations generated by $\opc_i$ are exactly those which avoid 1--3--2.
\end{lemma}
The proof is similar to that of Lemma \ref{lem:generate123}.

The three sets of operators are in a sense the same operator in different guises, as the next result shows.
\begin{theorem}\label{thm:operator1}
	Operators $\opa_i, \opb_i, \opc_i$ define bijections between the sets of Dyck paths of length $n$ with $k$ peaks, 1--2--3-avoiding permutations of $[n]$ with $k$ cycles, and 1--3--2-avoiding permutations of $[n]$ with $k$ cycles.
\end{theorem}
\begin{proof}\label{pf:operator1}
	It is clear that in each case, an object $o$ is uniquely defined by a sequence of numbers $i_1, \ldots, i_n$ such that $o$ is the object achieved by applying operators $i_1$ to $i_n$ in sequence to the empty object. 
	
	Furthermore, after applying operator $\opa_i, i>0$ on a path it will have $i$ excursions, so on that path we may apply operators $\opa_0$ to $\opa_i$. After applying $\opb_i, i>0$ to a permutation, the first $i-1$ cycles will be 1-cycles, so on that permutation we may apply operators $\opb_0$ to $\opb_i$. After applying $\opc_i, i>0$ to a permutation, there will be $i$ cycles whose smallest element is greater than all elements in all cycles to the right. On the resulting permutation we may apply operators $\opc_0$ to $\opc_i$.
	
	Finally, operators $\opa_0, \opb_0, \opc_0$ adds a peak, a cycle, and a cycle respectively, and no other operator adds peaks respectively cycles. Also, the number of excursions, the number cycles before the first non 1-cycle, and the number of cycles whose smallest element is greater than all elements in all cycles to the right, all increases by one.
	
	Therefore, if a sequence $i_1, \ldots, i_n$ defines a dyck path of length $n$ with $k$ peaks, it also defines an  1--2--3-avoiding permutations of $[n]$ with $k$ cycles, and an 1--3--2-avoiding permutations of $[n]$ with $k$ cycles.
\end{proof}

\subsubsection{Equivalence of avoidance of 2--1--3 and 3--1--2}
The ideas from the previous section is reused here, with some small alterations in the details. The proofs are very similar, and skipped. 
\subsubsection*{Operator set 2 for Dyck paths}
If a path $p$'s last decent is of length $k$, the operators $\opd_i$, $i=1,\ldots, k$ changes the postfix $S^k$ to $S^{i-1}NS^{k+1-i}$. Also define an operator $\opd_0$ by $\opd_0 p=p NS$. 

\subsubsection*{Operator set for 2--1--3}
Assume the first cycle in $\pi\in\mathcal{S}_n$ is $(c_1 c_2 \cdots c_k \cdots c_m)$, where $c_k$ is the first element that is not smaller than all elements $c_{k+1}, \ldots, c_m$ to the right of it in the cycle. Let $\ope_i$, $i=1,\ldots, k$, be the permutations whose cycles agree with those in $\pi$, except the first cycle which is $(c_1 c_2 \cdots c_i\ n+1\ c_{i+1} \cdots c_m)$. Also define $\ope_0$ by $\ope_0 \pi=(n+1)\pi$. 

\begin{lemma}
	The set of permutations generated by $\opc_i$ are exactly those who avoid 2--1--3.
\end{lemma}

\subsubsection*{Operator set for 3--1--2}
Assume the last cycle in $\pi\in\mathcal{S}_n$ is $(1 c_2 \cdots c_j \cdots c_m)$, where $c_j$ is the last element such that $c_{j+1}$ is greater. Let $\opf_i$, $i=1,\ldots, m+1-j$, be the permutations whose cycles agree with those in $\pi$, except the last cycle which is $(1 c_2 \cdots c_{m+1-i}\ n+1\ \cdots c_m)$. Also define $\opf_0$ by $\ope_0 \pi=(n+1)\pi$. 

\begin{lemma}
	The set of permutations generated by $\opf_i$ are exactly those who avoid 3--1--2.
\end{lemma}

\begin{theorem}\label{thm:operator2}
	Operators $\opd_i, \ope_i, \opf_i$ define bijections between the sets of Dyck paths of length $n$ with $k$ returns, 2--1--3-avoiding permutations of $[n]$ with $k$ cycles, and 3--1--2-avoiding permutations of $[n]$ with $k$ cycles.
\end{theorem}

From well known results about Dyck path statistics, it follow that the avoidance distributions for 2--1--3 and 3--1--2 (and 2--3--1, see below) are given by the Catalan triangle, and for 1--2--3 and 1--3--2 by the Narayana numbers, respectively. See Theorems \ref{theorem:1-3-2avoid} and \ref{theorem:2-1-3avoid} below.

\subsubsection{Equivalence of 2--3--1 and 3--1--2}
For a permutation $\pi$, define a mapping $\DDD$ by
\[\DDD(\pi)=\PAR^{-1}(\{\PAR(\pi)\}^{-1}).\]
\begin{theorem}
The mapping $\DDD(\pi)$ is an involution on the set of permutations in standard cycle form such that if $\pi$ has $k$ cyclic occurrences of 3--1--2 then $\DDD(\pi)$ has $k$ cyclic occurrences of 2--3--1. Furthermore, the number of cycles in $\pi$ and $\DDD(\pi)$ are the same. 
\end{theorem}
\begin{proof}
First we show that the number of cycles in $\pi$ and $\DDD(\pi)$ are the same. Elements in cycles are mapped to left-to-right minima in $\PAR(\pi)$ (and naturally vice versa). Also, left-to-right minima in a permutation in index form are mapped to left-to-right minima in the index form of the permutation's inverse.

Next we show that occurrence of 3--1--2 in a permutation (in index form) is mapped to occurrences of 2--3--1 in (the index form of) its inverse. Of course, a cyclic occurrence of a pattern in $\pi$ is in direct correspondence with an occurence in $\PAR(\pi)$.

Let $i<j<k$ be such that $\pi(i)>\pi(j)<\pi(k)$ and $\pi(i)>\pi(k)$, that is, $(i,j,k)$ is an occurrence of 3--1--2 in $\pi$. By rearranging we get $\pi(j)<\pi(k)<\pi(i)$ and $\pi^{-1}(\pi(j))<\pi^{-1}(\pi(k))>\pi^{-1}(\pi(i))$ and $\pi^{-1}(\pi(j))>\pi^{-1}(\pi(i))$, that is, $\pi(j),\pi(k),\pi(i)$ is an occurrences of 2--3--1 in $\pi^{-1}$.

That $\DDD$ is a involution is trivial.
\end{proof}
\begin{corollary}
	The number of 3--1--2-avoiding permutations of $[n]$ with $k$ cycles equals the  number of 2--3--1-avoiding permutations of $[n]$ with $k$ cycles: $f_{\mbox{\scriptsize 3--1--2}}(k,m,n)=f_{\mbox{\scriptsize 2--3--1}}(k,m,n)$.
\end{corollary}

\subsection{1--2--3}
The series starts 
\begin{align*}
F_{\mbox{\scriptsize1--2--3}}=&1+xz+(x+x^{2})z^{2}+\left((1+q)x+3x^{2}+x^{3}\right)z^{3}\\
&+\left((1+2q+2q^{2}+q^{4})x+(6+4q+q^{2})x^{2}+6x^{3}+x^{4}\right)z^{4}+O(z^{5}).
\end{align*}
For the avoidance distribution, see Theorem \ref{theorem:1-3-2avoid}.

\subsection{1--3--2}
The series starts 
\begin{align*}
F_{\mbox{\scriptsize 1--3--2}}(q,x,z)=&1+xz+(x+x^{2})z^{2}+\left((1+q)x+3x^{2}+x^{3}\right)z^{3}\\
&+\left((1+q+3q^{2}+q^{3})x+(6+4q+q^{2})x^{2}+6x^{3}+x^{4}\right)z^{4}+O(z^{5}).
\end{align*}

\begin{theorem}\label{theorem:1-3-2avoid}
The number of permutations of length $n$ with $k$ cycles that cyclically avoid 1--3--2 equals the number of Dyck paths of length $n$ with $k$ peaks, and the number of permutations of length $n$ with $k-1$ descents, and equals the Narayana numbers: 
\[f_{\mbox{\scriptsize 1--3--2}}(0,k,n)=f_{\mbox{\scriptsize1--2--3}}(0,k,n)=\frac{1}{k}{\binom{n-1}{k-1}}{\binom{n}{k-1}}.\]
\end{theorem}
\begin{remark}
That decents in 1--3--2-avoiding permutations \emph{also} are distributed as the Narayana numbers is known, see for instance \cite[Remark 2.5(b)]{Reif2003}.
\end{remark}

\subsection{2--1--3}
The series starts 
\begin{align*}
F_{\mbox{\scriptsize 2--1--3}}(q,x,z)=&1+xz+(x+x^{2})z^{2}+\left(2x+(2+q)x^{2}+x^{3}\right)z^{3}\\
&+\left((5+q)x+(5+2q+4q^{2})x^{2}+(3+2q+q^{3})x^{3}+x^{4}\right)z^{4}+O(z^{5}).
\end{align*}

From the Dyck path relation it follows that the avoidance distribution is given by the Catalan triangle.
\begin{theorem}\label{theorem:2-1-3avoid}
The number of permutations of length $n$ with $m$ cycles that avoids 3--1--2 equals the number of Dyck paths of length $n$ that have $m$ returns to the $x$-axis, i.e.~ the distribution is given by the Catalan triangle:
\[f_{\mbox{\scriptsize 3--1--2}}(0,m,n)=f_{\mbox{\scriptsize 2--3--1}}(0,m,n)=f_{\mbox{\scriptsize 2--1--3}}(0,m,n)=\frac{m}{2n-m}\binom{2n-m}{n}.\]
\end{theorem}

\subsection{2--3--1 and 3--1--2}\label{section:312and231}
The common series starts 
\begin{align*}
F_{\mbox{\scriptsize 3--1--2}}(q,x,z)=&F_{\mbox{\scriptsize 2--3--1}}(q,x,z)=1+xz+(x+x^{2})z^{2}+\left(2x+(2+q)x^{2}+x^{3}\right)z^{3}\\
&+\left((5+q)x+(5+3q+2q^{2}+q^{3})x^{2}+(3+q+2q^{2})x^{3}+x^{4}\right)z^{4}+O(z^{5}).
\end{align*}

The distribution of permutations avoiding one of these patters is given in Theorem \ref{theorem:2-1-3avoid}.

\subsection{3--2--1}
The series starts 
\begin{align*}
F_{\mbox{\scriptsize 3--2--1}}(q,x,z)=&1+xz+(x+x^{2})z^{2}+(2x+3x^{2}+qx^{3})z^{3}\\
&+\left((5+q)x+(9+2q)x^{2}+(3q+3q^{2})x^{3}+q^{4}x^{4}\right)z^{4}+O(z^{5}).
\end{align*}
Obviously a permutation that avoids 3--2--1 can have at most two cycles. From the basic bijection between cycles and permutations it follows from the classic result that the number of permutations with one cycle that avoids 3--2--1 are given by the Catalan numbers. Using the ideas use earlier in this section, the permutations with two cycles can be bijectively mapped to a set of Dyck paths.
\begin{theorem}\label{thm:321}
	The number of permutations with two cycles that avoid 3--2--1 and whose first cycle has $k$ elements is the same as the number of Dyck paths whose first peak is at height $k+1$. 
\end{theorem}
\begin{proof}[Sketch of proof]
	Starting with the path $NNSS$, Dyck paths with first peak at height at least two can be grown using the set of operators $\opd_i$. In the permutation domain, we start with $(2)(1)$. The growth operators insert $k+\tfrac{1}{2}$ last in the second cycle for all allowed values of $k$. If the last cycle has size one, also insert $n+\tfrac12$ last in the first cycle. In the path domain this corresponds to going from a path $N^nS^n$ to $N^{n+1}S^{n+1}$.  
\end{proof}
\begin{remark}
	The number of Dyck paths whose first peak is at height $k$ equals the number of Dyck paths with $k$ returns, and are thus given by the Catalan triangle. 
\end{remark}

\subsection{2--1--2', or ``2--1--3 or 3--1--2''}\label{section:2-1-2'}
The series starts
\begin{align*}
F_{\mbox{\scriptsize 2--1--2'}}=&1+xz+(x+x^{2})z^{2}+(2x+(1+2q)x^{2}+x^{3})z^{3}\\
&+\left((4+2q)x+(2+7q^{2}+2q^{3})x^{2}+(1+2q+q^{3}+2q^{4})x^{3}+x^{4}\right)z^{4}+O(z^{5})
\end{align*}
From Theorem \ref{theorem:nodash} it follows that
\begin{corollary}
\[G_{\mbox{\scriptsize2--1--2'}}(q,z)=F_{\mbox{\scriptsize2--1--2'}}(q,1,z).\]
\end{corollary}

It is easy to see, \cite{Kita2005}, that $g_{\mbox{\scriptsize2--1--2'}}(0,n)=2^{n-1}$. This is extended here.
\begin{theorem}
\[f_{\mbox{\scriptsize2--1--2'}}(0,m,n)=2^{n-m-1}\mbox{ for }1\leq m\leq n-1\mbox{ and }f_{\mbox{\scriptsize2--1--2'}}(0,n,n)=1.\]
\end{theorem}
\begin{proof}
It is trivial to check that the statement holds for short permutations.
Assume the statement is true for $n$-permutations. From $\mathcal S_{n}$ we can get $\mathcal S_{n+1}$ in the following way:

For each permutation $\pi$ insert the element $n+1$ in the same cycle as and after $1, 2,  \ldots, n$. Additionally we can add the fixed point $(n+1)$. Now consider a 2--1--2'-avoiding permutation $\pi$ with $m$ cycles, $1\leq m\leq n-1$. 

There are exactly two places to insert $n+1$ into $\pi$ to give a 2--1--2' avoiding permutation; before or after the element $n$. 
(As $\pi$ is 2--1--2'-avoiding $n$ cannot be (assuming $n\geq3$) a fixed point so it is possible to insert $n+1$ before $n$.) 

Furthermore, it is easy to see that there is exactly one 2--1--2'-avoiding permutation of length $n$ and $n-1$ cycles, and that the sole permutation with $n$ cycles is 2--1--2'-avoiding.
\end{proof}

The distribution of one occurrence of 2--1--2' is obtained in the same way.
\begin{theorem}~
\begin{enumerate}[\it i.]
\item $f_{\mbox{\scriptsize2--1--2'}}(1,m,n)=2^{n-m-2}\mbox{ for }1\leq m\leq n-3$. 
\item $f_{\mbox{\scriptsize2--1--2'}}(1,n-1,n)=1$. 
\item $f_{\mbox{\scriptsize2--1--2'}}(1,n,n)= f_{\mbox{\scriptsize2--1--2'}}(1,n-2,n)=0$.
\end{enumerate}
\end{theorem}
\begin{proof}
We will use induction to show \emph i. Let $\mathcal P^{m}_{n}$ denote the set of permutations with $m$ cycles and one cyclic occurrence of 2--1--2'. 
 To start the induction, the two permutations in $\mathcal P^{1}_{4}$ are (1324) and (1324). 

From $\mathcal P^{1}_{n}$ we get $\mathcal P^{1}_{n+1}$ in following way. For each $\pi\in \mathcal P^{1}_{n}$ add 1 to every element except 1, and insert 2 \emph{a}) after 1 or \emph{b}) last in the cycle. It is easy to see that the resulting permutations are in $\mathcal P^{1}_{n+1}$, and that inserting 2 in any other place gives a permutation with more than one occurrence of 2--1--2'.

From $\mathcal P^{m}_{n}$ we get $\mathcal P^{m+1}_{n+1}$ in following way. For each $\pi\in \mathcal P^{m}_{n}$ add 1 to every element and add the cycle (1). It is easy to see that the resulting permutations are in $\mathcal P^{m}_{n+1}$, and that inserting 1 in any other place give a permutation with more than one occurrence of 2--1--2'.

For \emph{ii}, it is easily verified that the two permutations in $\mathcal P^{n-1}_{n}$ are 
\[(n-1)(n-2\ n)(n-3)\cdots (2)(1) \mbox{ and } (n)(n-2\ n-1)(n-3)\cdots (2)(1).\]

For \emph{iii}, consider one of the above permutations. There is no place to insert the element $n+1$ into any cycle without adding occurrences of 2--1--2'. Finally, $(n)(n-1)\cdots (2)(1)$ avoids 2--1--2'.
\end{proof}
In the same way the distribution of permutations with two occurrences of 2--1--2' may be found.
\begin{theorem}~
\begin{enumerate}[\it i.]
\item $f_{\mbox{\scriptsize2--1--2'}}(2,m,n)=2^{n-m-1}\mbox{ for } 1\leq m \leq n-4$. 
\item $f_{\mbox{\scriptsize2--1--2'}}(2,n-2,n)=7$ for $n\geq  4$. 
\item $f_{\mbox{\scriptsize2--1--2'}}(2,n-1,n)=1$ for $n\geq  4$. 
\item $f_{\mbox{\scriptsize2--1--2'}}(2,k,n)=0$ otherwise.
\end{enumerate}
\end{theorem}
Also, we conjecture that in general it holds that 
\[g_{\mbox{\scriptsize2--1--2'}}(k,n)=2g_{\mbox{\scriptsize2--1--2'}}(k, n-1)\mbox{ for } 0\leq k \leq n-2.\]

\subsection{1--2--1', or ``2--3--1 or 1--3--2''}\label{section:1-2-1'}
The series starts
\begin{align*}
F_{\mbox{\scriptsize1--2--1'}}=&1+xz+(x+x^{2})z^{2}+\left((1+q)x+(2+q)x^{2}+x^{3}\right)z^{3}\\
&+\left((1+q+2q^{2}+2q^{3})x+(3+2q+4q^{2}+2q^{3})x^{2}+(3+q+2q^{2})x^{3}+x^{4}\right)z^{4}+O(z^{5}).
\end{align*}
From Theorem \ref{theorem:nodash} it follows that
\begin{corollary}
\[G_{\mbox{\scriptsize1--2--1'}}(q,z)=F_{\mbox{\scriptsize1--2--1'}}(q,1,z).\]
\end{corollary}

Some distributional results are next.
\begin{theorem}~
\begin{enumerate}[\it i.]
\item $f_{\mbox{\scriptsize1--2--1'}}(0,k,n)=\binom{n-1}{k-1}$ for $1\leq k \leq n$.
\item $f_{\mbox{\scriptsize1--2--1'}}(1,k,n)=\binom{n-2}{k-1}$ for $1\leq k \leq n-1$, $n\geq 3$.
\item $f_{\mbox{\scriptsize1--2--1'}}(2,k,n)=2\binom{n-2}{k-1}$ for $1\leq k \leq n-1$, $n\geq 4$.
\end{enumerate}
\end{theorem}
\begin{proof}
For the first assertion, we can arrange the 1--2--1'-avoiding permutations into a Pascal-like triangle:\\

\begin{center}
\setlength{\tabcolsep}{-2.0mm}
\setlength{\arraycolsep}{0mm}
\begin{tabular}{p{2.0cm}p{2.0cm}p{2.0cm}p{2.0cm}p{2.0cm}p{2.0cm}p{2.0cm}p{2.0cm}p{3cm}}
& & & & 
\centering (1) 
& & & & \\

& & &
\centering (12) & &
\centering (2)(1) 
& & & \\

& &
\centering (123) & &
\centering
\begin{tabular}{c}
(3)(12)\\
(2)(13) 
\end{tabular}& &
\centering (3)(2)(1) 
& & \\

&
\centering (1234) & &
\centering 
\begin{tabular}{c}
(4)(123)\\ 
(3)(124)\\
(2)(134) 
\end{tabular}& &
\centering 
\begin{tabular}{c}
(4)(3)(12)\\
(3)(2)(14)\\
(4)(2)(13)
\end{tabular}& &
\centering (4)(3)(2)(1) &  \\

\centering (12345) & &
\centering 
\begin{tabular}{c} 
(5)(1234)\\ 
(4)(1235)\\ 
(3)(1245)\\
(2)(1345) 
\end{tabular}& &
\centering 
\begin{tabular}{c}
(5)(4)(123)\\ 
(5)(3)(124)\\
(5)(2)(134)\\
(3)(2)(145)\\
(4)(2)(135)\\
(4)(3)(125) 
\end{tabular}& &
\centering 
\begin{tabular}{c} 
(4)(3)(2)(15)\\
(5)(3)(2)(14)\\
(5)(3)(2)(13)\\
(5)(4)(2)(12) 
\end{tabular}& &
\centering 
(5)(4)(3)(2)(1)
\end{tabular}
\end{center}
Each set of permutations, say of length $n$ and with $k$ cycles, are the union of 
\begin{enumerate}[\it i.]
\item the set achieved from the set above to the left with the added cycle $(n+1)$, and 
\item the set achieved from the set above to the right with the element $n+1$ added last in the cycle including 1.
\end{enumerate}
We conclude that the sizes of the sets in the triangle are the binomial coefficients.

The same idea works for permutations with one and two cyclic occurrence of 1--2--1'. In first case start at level three with (132) and (2)(13), in the second at level four with (1243) and (1423); (2)(143), (24)(13), (3)(142) and (34)(12); and (34)(2)(1) and (3)(24)(1).
\end{proof}

\begin{corollary}
For $n\geq 2$, $g_{\mbox{\scriptsize 1--2--1'}}(1,n)=2^{n-2}$. For $n\geq 3$, $g_{\mbox{\scriptsize 1--2--1'}}(2,n)=2^{n-1}$.
\end{corollary}

We also conjecture that the following results hold.
\begin{enumerate}[\it i.]
\item $f_{\mbox{\scriptsize1--2--1'}}(3,k,n)=2\binom{n-2}{k-1}+2\binom{n-3}{k-1}$ for $1\leq k \leq n-1$, $n\geq 5$.
\item $f_{\mbox{\scriptsize1--2--1'}}(4,k,n)=3\binom{n-2}{k-1}+2\binom{n-3}{k-1}$ for $1\leq k \leq n-1$, $n\geq 6$.
\item $f_{\mbox{\scriptsize1--2--1'}}(5,k,n)=2\binom{n-2}{k-1}+7\binom{n-3}{k-1}$ for $1\leq k \leq n-1$, $n\geq 7$.
\item $f_{\mbox{\scriptsize1--2--1'}}(6,k,n)=4\binom{n-2}{k-1}+8\binom{n-3}{k-1}+4\binom{n-4}{k-1}$ for $1\leq k \leq n-1$, $n\geq 8$.
\item $f_{\mbox{\scriptsize1--2--1'}}(7,k,n)=2\binom{n-2}{k-1}+12\binom{n-3}{k-1}+6\binom{n-4}{k-1}$ for $1\leq k \leq n-1$, $n\geq 9$.
\end{enumerate}


\bibliographystyle{abbrv}

\end{document}